\documentclass[oneside,12pt]{amsart}
\usepackage{multirow} 
\usepackage{amscd,amsmath,latexsym,bm,pdfpages,verbatim,amsthm}
\usepackage[mathcal]{euscript}
\usepackage{graphicx}
\usepackage{amssymb}          
\usepackage{epstopdf}
\usepackage[all]{xy}

\newtheorem{thm}{Theorem}[section]
\newtheorem{cor}[thm]{Corollary}

\newtheorem{prop}[thm]{Proposition}
\numberwithin{equation}{section}
\theoremstyle{definition}
\newtheorem{defin}[thm]{Definition}
\theoremstyle{definition}

\newtheorem{remark}[thm]{Remark}
\theoremstyle{conjecture}

\theoremstyle{remark}
\newtheorem*{rem}{Remark}
\theoremstyle{remark}
\newtheorem*{remex}{Example}

\def\nd{\noindent}

\def\skp#1{\vskip#1cm\relax}

\begin{document}

\title[On the rational type of moment-angle complexes]{On the rational type of moment-angle complexes}

\author[A.~Bahri]{A.~Bahri}
\address{Department of Mathematics,
Rider University, Lawrenceville, NJ 08648, U.S.A.}
\email{bahri@rider.edu}

\author[M.~Bendersky]{M.~Bendersky}
\address{Department of Mathematics
CUNY,  East 695 Park Avenue New York, NY 10065, U.S.A.}
\email{mbenders@xena.hunter.cuny.edu}

\author[F.~R.~Cohen]{F.~R.~Cohen }
\address{Department of Mathematics,
University of Rochester, Rochester, NY 14625, U.S.A.}
\email{cohf@math.rochester.edu}

\author[S.~Gitler]{S.~Gitler}
\address{Department of Mathematics,
Cinvestav, San Pedro Zacatenco, Mexico, D.F. CP 07360 Apartado
Postal 14-740, Mexico. Colegio Nacional. Institute for Advanced Studies, Princeton.} 
\email{Deceased}



 \thanks{The first author was supported in part by a Rider University Summer Research Fellowship and 
grant number 210386 from the Simons Foundation, the third was partially supported by the Institute for Mathematics and its Applications and the fourth author received support from {\sc CONACYT}, Mexico.}

\begin{abstract}
In this note, it is shown that the only moment-angle complexes  which are rationally elliptic are those which are products of odd spheres and a disk.
\skp{0.2}
\end{abstract}

\maketitle

\section{Introduction}

\noindent F\'elix and Halperin showed,  \cite{fh} and \cite{fht}, that there is a dichotomy for simply-connected 
finite $CW$-complexes $X$. Their theorem is the following.
\begin{thm}Either
\begin{enumerate}
\item[(1)] {\bf  $\pi_*(X) \otimes  \mathbb{Q}$} is a finite $\mathbb{Q}$-vector space, in which case $X$ is called
rationally elliptic or,
\item[(2)]  {\bf  $\pi_*(X) \otimes  \mathbb{Q}$} grows exponentially, in which case $X$ is called  rationally hyperbolic.
\end{enumerate}
\end{thm}

\noindent  Next,  recall the definition of the {\em moment-angle complex} $Z(K;(D^2;S^1))$ from \cite{bbcg1}, \cite{5}.\\
\begin{defin} Let $(D^2, S^1)$ be the pair of a 2-disk and its boundary circle, and $K$ be a finite abstract simplicial complex with $n$ vertices. Then
$Z(K;(D^2;S^1))$ is a subspace of $(D^2)^n$ defined as the  union over all simplices $\sigma \in K$ of subspaces of $(D^2)^n$
$$D(\sigma)=\{(x_1, ..., x_n) \mid x_i \in S^1 \mbox{  if  } i \not\in \sigma \},$$
$Z(K;(D^2;S^1))$ is a 2-connected finite $CW$-complex.
\end{defin}

\nd In this note we obtain:

\begin{thm}\label{thm:rational.homotopy.type}
The only moment-angle complexes $Z(K;(D^2;S^1))$ which are rationally elliptic, are those which are a product of odd spheres and a disk.
\end{thm}

\begin{rem} After an inquiry from one of the authors, A.~Berglund,  \cite{berg}, supplied an alternative proof of Theorem \ref{thm:rational.homotopy.type}. Also, the authors have
learned of the related work of  Gery Debongnie \cite{deb}, done in the context of the complement $M(A)$ of a subspace arrangement $A$. The moment-angle complex $Z(K;(D^2,S^1))$ is homotopy equivalent to the complement of a subspace arrangement given by coordinate planes. For arrangements with a geometric lattice with subspaces of codimension at least two, Debongnie classified those M(A) which are rationally elliptic as homotopy equivalent to a finite product of odd spheres. Though Debongnie's theorem is more general, 
the purpose of the current paper is to give a proof for moment-angle complexes from first principles.
 In his thesis \cite{gur}, Michael Gurvich determined that the toric manifolds which are rationally elliptic
arise from simple polytopes which are products of simplices. A short homotopy argument proves Theorem
\ref{thm:rational.homotopy.type} for $K$ dual to the boundary of such simple polytopes.
\end{rem}

\section{Minimal non-faces and abstract simplicial complexes} 
Let $[n]$ denote an abstract set of vertices $\{v_1, ..., v_n\}$. The next definition yields a description of an abstract simplicial complex in term of its ``missing'' faces.

\begin{defin}\label{defin:def_M}

A family $M=\{m_1, ..., m_k\}$ of subsets of $[n]$ satisfying:
\begin{enumerate}
\item[(1)] $|m_i| > 1$ and,
\item[(2)] $m_i \not\subset m_j$ for any $i \not= j \in \{1,2,..   .,k\}$ 
\end{enumerate}
is called a set of minimal non-faces. 
\end{defin}

\noindent Let $\nu  = \bigcup_{i=1}^{n}m_i$.      
 Associated to $M$ are two abstract simplicial complexes.
$$K(M,[n])=\{\sigma \subset [n] \mid m_i \not\subseteq \sigma \mbox{ for all  } i=1,...,k\}$$
$$K(M,\nu )=\{\sigma \subset \nu \mid m_i \not\subseteq \sigma \mbox{ for all  } i=1,...,k\}$$
which agree if $\nu=[n]$. If $|\nu| < n$, $K(M, \nu)$, is called the {\em reduced\/} simplicial complex corresponding to $M$. The empty set $\varnothing$ is 
considered to be in both simplicial complexes.\\

\noindent   In a simplicial complex $K$, 
a minimal non-face is a sequence  of  vertices $Q=(v_{i_{1}},v_{i_{2}},\ldots, v_{i_{q}})$
so that $Q \not\in K$, but every proper subsequence of $Q$ is a simplex of $K$.

\begin{remark}

If $K$  is  an abstract simplicial complex with $n$ vertices and $M$ is its set of minimal non-faces, then there is a homeomorphism of underlying simplicial complexes
$$K \longrightarrow  K(M,[n]),$$
\end{remark}

\noindent Recall next that the join of two disjoint simplicial complexes $K_1$ and $K_2$, denoted by $K_1 \ast K_2$, is defined by
$$K_1\ast K_2=\{ \sigma_1 \cup  \sigma_2 \mid \sigma _1 \in  K_1, \sigma_2 \in K_2\}.$$

\begin{prop}\label{prop:disc}
Let $|\nu | < n$ and set $n'=|\nu|$, then there is a simplicial isomorphism:
$$K(M,[n]) \longrightarrow K(M, \nu) \ast \Delta^{n-n'-1}$$
where $\Delta^{n-n'-1}$ is the simplex with $n-n'$ vertices $ \{v_{i_1}, ..., v_{i_{n-n'}}\} = [n]-\nu$.
\end{prop}

\begin{proof}
The sets $\nu$ and $[n]-\nu$ are disjoint, so a simplicial isomorphism is given by
$$\sigma \longrightarrow \sigma_1 \cup \sigma_2$$

\nd  where $\sigma_1=\sigma \cap K(M,\nu)$, $\sigma_2=\sigma \cap K(M,[n]-\nu)$.
Notice that $\Delta^{n-n'-1}$ is a simplex because every subset of it  is a simplex of $K(M,[n])$. 
\end{proof}

\begin{defin}\label{defin:G(M)}

Given $M$ as in Defiition \ref{defin:def_M}, a graph $G(M)$ is defined with its vertices the $m_i$ and an edge joining $m_i$ and $m_j$ if $m_i \cap m_j \neq \varnothing$.\\

\noindent Let $G(M)=\{C_1(M), \ldots , C_l(M)\}$ be the connected components of $G(M)$ and let
$M_i \subset M$ be the set of $m_j$ in $C_i(M)$.
\end{defin}

\begin{prop}
There is a simplicial isomorphism
$$K(M, \nu) \xrightarrow{\varphi} K(M_1, \nu_{1}) \ast K(M_2,\nu_2) 
\ast \cdots \ast K(M_l,\nu_l).$$
\end{prop}
\begin{proof}

Let $\sigma \in K(M, \nu)$ and $\sigma_i$ be the part of $\sigma$ which lies in $K(M_i,\nu_i)$; notice that
$\sigma_i$ is a simplex in $K_i$. Set
$$\phi(\sigma) = \sigma_1 \cup  \sigma_2 \cup \ldots \cup  \sigma_l$$
\noindent where $\sigma_1 \cup  \sigma_2 \cup \ldots \cup  \sigma_l \in K(M_1,\nu_1) \ast \cdots 
\ast K(M_l,\nu_l)$.
Conversely, every such union is a simplex of $K(M,\nu)$ since it does not contain any of the $m_i \in M$. \end{proof}

\begin{prop}
Let $m \in M$ be a minimal non-face, Then $K(M, |m|)$ is isomorphic to the boundary of a simplex, $\partial \Delta (|m|-1)$.
\end{prop}

\begin{proof} Every proper subsequence of $m$ is a simplex in $K(M, |m|)$, and then isomorphic to the boundary 
$\partial \Delta (|m| -1)$. \end{proof}

\begin{cor}\label{cor:elliptic}
If the minimal non-faces of $M$ are pairwise disjoint, then there is a simplicial isomorphism 
$$K(M,\nu) \cong \partial \Delta (|m_1|-1)\ast \ldots \ast \partial \Delta (|m_k|-1)$$
where $|M|=k$.
\end{cor}

\section{The dichotomy for moment-angle complexes}

\noindent The following properties of moment-angle complexes may be found in \cite{bbcg1} and \cite{5}.
\begin{align}
Z(K_1 \ast K_2; (D^2, S^1)) &\cong  Z(K_1;(D^2, S^1)) \times Z(K_2;(D^2, S^1))\\
Z(\Delta^k; (D^2, S^1)) &\cong  D^{2k+2}\\
Z(\partial \Delta^k; (D^2, S^1)) &\cong  S^{2k+1}
\end{align}
\skp{0.1}

\nd  Proposition \ref{prop:disc} and Corollary \ref{cor:elliptic} imply now that if all the minimal non-faces
of $K$ are pairwise disjoint, $Z(K;(D^{2},S^{1}))$  is the product of odd spheres and a disc and hence is rationally
elliptic

 The case not covered by Corollary 2.7 is addressed next.

\begin{defin}
Let $\mathcal{A}_{m}$ be the collection of all simplicial complexes on $m$ vertices which have
a pair of intersecting minimal non-faces, but no proper full subcomplex with that property.
\end{defin}

\begin{remex}
Let $m=4$ and $K$ have minimal non-faces corresponding to relations in the Stanley-Reisner ring:
$v_{1}v_{2}v_{3}$, $v_{1}v_{2}v_{4}$ and $v_{1}v_{4}$. Here, $K$ has no proper full subcomplex
with intersecting non-faces.
\end{remex}

\begin{prop}\label{prop:am}
Let $K \in \mathcal{A}_{m}$, then $Z(K;(D^{2},S^{1}))$ has a wedge of odd spheres as a retract.
\end{prop}
\begin{proof}
Suppose that $K$ has minimal intersecting non-faces corresponding to the following relations in the 
Stanley-Reisner ring
$$v_{1}\cdots v_{k}w_{1}\cdots w_{t}\quad\text{and}\quad u_{1}\cdots u_{r}w_{1}\cdots w_{t}.$$

\nd (Notice that minimality dictates that $k$, $t$ and $r$ are all $\geq 1$.) 
It follows that the vertex set of $K$ must be
\begin{equation}\label{eqn:vertices}
\big\{v_{1},\ldots, v_{k},u_{1},\ldots,u_{r},w_{1},\ldots,w_{t}\big\}
\end{equation}

\nd for otherwise, removing a vertex from $K$, which is not among these, will produce a proper full
subcomplex contradicting $K \in \mathcal{A}_{m}$.  Next, setting
$$I =\big\{v_{1},\ldots, v_{k},w_{1},\ldots,w_{t}\big\} \quad\text{and}\quad 
J = \big\{u_{1},\ldots,u_{r},w_{1},\ldots,w_{t}\big\}$$

\nd gives retractions off $Z(K;(D^{2},S^{1}))$:
$$Z_{K_{I}} = S^{2(k+t) -1} \quad\text{and}\quad Z_{K_{J}} = S^{2(r+t) -1}.$$

\nd corresponding to the full subcomplexes $K_{I}$ and $K_{J}$,  \cite[Theorem $2.2.3$]{denham.suciu} . The stable splitting theorem of \cite{bbcg1} distinguishes these two spheres. This gives a map
$$S^{2(k+t) -1} \vee S^{2(r+t) -1} \longrightarrow Z(K;(D^{2},S^{1})).$$

\nd It remains to show that rationally, no cells are attached to this wedge of spheres inside $Z(K;(D^{2},S^{1}))$.
Now, the results of \cite{bbcg1} imply that all non-trivial attaching maps to this wedge of spheres must be 
{\em stably trivial\/}.
The Hilton-Milnor theorem, \cite[Theorem $4.3.2$]{jn}, gives
\begin{align*}\pi_{n}(S^{2(k+t) -1} \vee S^{2(r+t) -1})&\cong \pi_{n}(S^{2(k+t) -1})\oplus  \pi_{n}(S^{2(r+t) -1})\oplus 
\pi_{n}\big(\Sigma(S^{2(k+t) -2} \wedge S^{2(r+t) -2})\big)\\
&\oplus_{j \geq 2} \pi_{n}\big(\Sigma(S^{2j(k+t) -j} \wedge S^{2(r+t) -1})\big).                       
\end{align*}

\nd The rational homotopy groups of spheres is well known. The only stably trivial non-trivial classes occur in the groups
$\pi_{4q-1}(S^{2q})$.  In the decomposition above, this requires
$$n \geq  4(2k+3t+r-1)-1.$$

\nd The vertex set of $K$ is given by \eqref{eqn:vertices} and so the largest cell possible in $Z(K;(D^{2},S^{1}))$
has dimension $2(k+r+t)-1$. Now
$$2(k+r+t)-1 < 4(2k+3t+r-1)-1$$

\nd because $k$, $t$ and $r$ are all $\geq 1$. So rationally, no non-trivial attaching map is possible. \end{proof}

An induction argument now gives the result.

\begin{thm}\label{thm:induction}
Let $K$ be a simplicial complex which contains a pair of minimal intersecting non-faces, then $Z(K;(D^{2},S^{1}))$
is rationally hyperbolic.
\end{thm}

\begin{proof}
It is straightforward to check that all simplicial complexes on three vertices, which have pairwise intersecting 
non-faces have a wedge of spheres as a retract and so are rationally hyperbolic.

Suppose by way of induction, that all simplicial complexes with fewer than $m$ vertices,
which have pairwise intersecting non-faces, have a wedge of spheres as a rational retract.  Let
$K$ be a simplicial complex on m vertices, which has pairwise intersecting non-faces. If
$K \in \mathcal{A}_{m}$, the result is true for $K$ by Proposition \ref{prop:am}. If $K \notin \mathcal{A}_{m}$,
then $K$ has a proper full subcomplex $L$ which has a pair of intersecting non-faces.  Finally, the induction
hypothesis and \cite[Theorem $2.2.3$]{denham.suciu} now imply the result. \end{proof}

\bibliographystyle{amsalpha}

\end{document}